\documentclass[11pt]{article}
\usepackage[letterpaper,margin=1in]{geometry}
\usepackage[T1]{fontenc}
\usepackage{lmodern,amsmath,amssymb,amsthm,mathtools,microtype,booktabs}
\usepackage{xcolor}
\usepackage[colorlinks=true,allcolors=blue!45!black]{hyperref}
\newtheorem{theorem}{Theorem}
\newtheorem{corollary}[theorem]{Corollary}
\newtheorem{lemma}[theorem]{Lemma}
\newtheorem{proposition}[theorem]{Proposition}
\newcommand{\ind}{\mathbf1}
\newcommand{\qb}[2]{\genfrac{[}{]}{0pt}{}{#1}{#2}_q}
\newcommand{\mass}[1]{\lvert #1\rvert}
\newcommand{\arx}[2]{\emph{#2}. \href{https://arxiv.org/abs/#1}{arXiv:#1}}
\title{A sharp bounded bijection for tight cylindric partitions}
\author{Scott D. Hughes}
\date{}
\hypersetup{pdftitle={A sharp bounded bijection for tight cylindric partitions},
  pdfauthor={Scott D. Hughes}}
\begin{document}
\maketitle
\vspace{-2.5em}
\begin{abstract}
Fix a cyclic profile of coprime rank \(r\) and level \(\ell\), both at least
three, and an entry bound \(n\ge1\). Let \(d\) be the deficit from the
maximum volume of a tight cylindric partition. For \(d\le n\min(r,\ell)\),
the one-row combinatorial \(R\) matches the tight cylinders whose last
slice lies inside the maximizer's last slice with the partitions of \(d\)
in which each part occurs fewer than \(r\) times. The inverse preserves
the profile and entry bound. We identify every partition missed at the
next deficit. A prefix test decides whether any cylinder lies outside this
sector and, if so, identifies the unique one of least deficit.
This explains the leading coefficient gap in Warnaar's \(Q_{n,c}\).
We also determine the first correction to the vertex terms of the
published bounded determinant.
\end{abstract}

\section{The bounded bijection and the exceptional chain}
An outgoing profile \(c=(c_0,\ldots,c_{r-1})\) consists of nonnegative
integers of sum \(\ell\). At height \(h\), let \(x_i^h\) be the number of entries
at least \(h\) in row \(i\) of a cylindric partition. Equivalently, the
objects are arrays of nonnegative integers with
\[
 x_i^{n+1}=0,\qquad x_i^{h+1}\le x_i^h,\qquad
 x_{i+1}^h\le x_i^h+c_i\quad(1\le h\le n),
\]
with row indices modulo \(r\). Their volume is \(W=\sum_{h,i}x_i^h\).
A cylinder is \emph{tight} if at every height \(h\) some row has
\(x_i^h=x_i^{h+1}\). Thus, for each part size, at least one row omits that part
\cite[Definition 3]{KR}. Cyclically rotating both rows and profile
preserves the objects and their volumes.

Set \(s_j=\sum_{i<j}c_i\). Throughout the paper we use the unique rotation
with
\[
 \delta_j=rs_j-\ell j>0\quad(0<j<r),\qquad \delta_0=\delta_r=0.
 \tag{1}\label{scores}
\]
Indeed, the proper partial sums \(rs_j-\ell j\), indexed by \(0\le j<r\),
are distinct modulo \(r\) since \(\gcd(r,\ell)=1\). Start at their minimum.
Let \(x^*\) be the unique tight cylinder of greatest volume \(M\), computed
in Section~\ref{geometry}. For a tight cylinder \(x\) its \emph{deficit}
is \(d=M-W(x)\). It belongs to the \emph{ordinary sector} if its last slice
is contained in that of \(x^*\), that is, \(x_i^n\le x_i^{*,n}\) for every
\(i\). Otherwise we call the cylinder \emph{signed}.

\begin{theorem}\label{main}
Let \(r,\ell\ge3\) be coprime, \(n\ge1\), and \(m=\min(r,\ell)\).
\textup{(i) Bounded inverse.} For each \(0\le d\le mn\) there is an
explicit bijection from the ordinary sector at deficit \(d\) to the
ordinary partitions of \(d\) with multiplicities less than \(r\). Each
inverse image has entries at most \(n\). At \(d=mn+1\) the map is still
injective and misses the following partition for every profile:
\[
 \begin{cases}
 (n+1,\underbrace{n,\ldots,n}_{r-1}),&r<\ell,\\
 (\ell n+1),&\ell<r .
 \end{cases}
\]
This is the sole missing partition except when \(n=1\), \(\ell=r+1\)
and the last entry of the canonical profile is one. In that case the
single part \((r+1)\) is also missed.

\textup{(ii) First signed cylinder.} Put \(E=\{j:1\le j<r,\ s_j<\ell\}\).
If \(E\) is empty, equivalently
\(c=(\ell,0,\ldots,0)\) in the canonical rotation, no signed cylinder
exists. Otherwise let \(j_*\) minimize \(\delta_j\) on \(E\), and put
\(\kappa=\delta_{j_*}\). The first signed deficit is \(n\kappa\).
Exactly one cylinder attains it. In the rotating order of \eqref{frame}
every slice has difference \(v_{j_*}\) from the maximum, where
\[
 v_j=(\underbrace{s_j,\ldots,s_j}_{r-j},
       \underbrace{s_j-\ell,\ldots,s_j-\ell}_{j}).
 \tag{2}\label{exception}
\]
For \(n\ge2\), raising the first lower-block entry of slice one by one
gives the unique signed cylinder at \(n\kappa+1\).
In particular, for the unique \(1\le j<r\), \(t\in\mathbb Z\) satisfying
\(rt-\ell j=1\), a signed cylinder occurs at \(d=n\) exactly when \(s_j=t\).
\end{theorem}

For the coefficient consequence, let \(C_n(q)\) enumerate all bounded cylinders
by volume, and put \((q;q)_n=\prod_{h=1}^n(1-q^h)\). Set \(Q_{0,c}=1\) and
define \(Q_{n,c}\) by triangular inversion of
\[
 F_n(q)=(q;q)_n C_n(q)=\sum_{k=0}^n\qb nk Q_{k,c}(q),\qquad
 \qb nk=\frac{(q;q)_n}{(q;q)_k(q;q)_{n-k}} .
 \tag{3}\label{tri}
\]
These are the coprime \(Q_{n,c}\) of \cite[(2.4)]{W}. Indeed, if
\(G_c(u,q)\) marks the largest entry, then \(\sum_{n\ge0}C_n(q)u^n
=G_c(u,q)/(1-u)\). Multiplying by \((u;q)_\infty\) and using the
\(q\)-binomial theorem yields \eqref{tri}, since the multiplier on \(G_c\)
is then \((uq;q)_\infty\). Outgoing offsets cyclically permute the profile
coordinates in \cite{W}.

\begin{corollary}\label{coeff}
Put \(L=\sum_{j=0}^{r-1}s_j\),
\(A=(r-1)(\ell-1)/2\), \(B=L-(r-1)(\ell+1)/2\), and \(d_n=An^2+Bn\).
The degree of \(Q_{n,c}\) at infinity is \(d_n\), and
\[
 z^{d_n} Q_{n,c}(1/z)\equiv
   \sum_{j=0}^{r-1} z^{n\delta_j}+z^{n+1}\pmod{z^{n+2}}.
 \tag{4}\label{formula}
\]
For \(n\ge2\), exactly \(n-1\) consecutive coefficients below the top
vanish if \(s_j=t\) in Theorem~\ref{main}, and exactly \(n\) otherwise.
The next coefficient is one. For \(n=1\), the coefficients at deficits
one and two are \(\ind_{\{\exists j:\delta_j=1\}}\) and
\(1+\ind_{\{\exists j:\delta_j=2\}}\), respectively.
\end{corollary}

We work with expansions of rational functions at infinity.
Warnaar credits Welsh with polynomiality and the value at \(1\)
\cite[following Conjecture 8.4]{W}. Removing and reinserting common copies
of each part gives the bounded tightening bijection of \cite[\S5]{KOS}.
Reading a tight level-one cylinder by columns gives the multiplicity model
\cite[Proposition 17]{KR}. Highest lifts and reconstruction from paths
appear in \cite[Propositions 2.11 and 2.13]{FLOTW} and
\cite[\S4 and Appendix C]{FW}. Tingley's tight-abacus bijection gives paths
with an eventual ground-state tail \cite[Theorem 5.1]{T}.
Huang, Lau and Ono give a bounded
cycle correspondence for balanced profiles \cite[Proposition~8.1, (8.1)]{HLO}.
The row exchange and its inverse on padded words are the published
combinatorial \(R\) \cite[\S2.3]{FOY}. Eventual stabilization with enough
vacuum sites is proved in \cite[Lemma 2.7]{FOY}. Here we determine exactly
how far the inverse fits the original entry bound and terminal profile cap,
and identify every obstruction at the next deficit. The signed-cylinder
classification then explains the coefficient gap. The bounded determinant
already implies \eqref{formula}. At rank three
and level four there are also explicit bounded sums \cite[Theorem 1.2]{CW}.

\section{The maximum and its terminal boundary}\label{geometry}
In slice \(h\), set \(a_i^h=c_i+x_i^h-x_{i+1}^h\). The slack vector
ranges over the nonnegative integer vectors of sum \(\ell\). Write
\(x_i^h=w_i^h+\tau_h\), with
\(w_i^h=s_i-\sum_{k<i}a_k^h\). For fixed slacks, tightness fixes
\(\tau_n=-\min_i w_i^n\) and
\(\tau_h=\tau_{h+1}+J(w^h,w^{h+1})\), where \(J(w,w')=\max_i(w_i'-w_i)\).
Hence
\[
 W=\sum_{h,i} w_i^h+r\sum_{h<n}h J(w^h,w^{h+1})-rn\min_i w_i^n .
\]
This is convex on the product of real slack simplices. At a vertex
\(\ell e_k\), the heights are \(s_i-\ell\ind_{i>k}\) and their minimum is
\(s_{k+1}-\ell\), using \(s_r=\ell\). The jump from vertex \(k\) to \(l\)
is \(\ell\ind_{k<l}\). For paths of \(n\) vertices ending at \(k\), the
maximum of \(\sum_h k_h+r\sum_{h<n}h\ind_{k_h<k_{h+1}}\) is
\((r-1)\binom n2+nk\). This holds for \(n=1\). In the induction the previous
endpoint \(l\) maximizes \(l+r\ind_{l<k}\), uniquely at \(k-1\bmod r\).
Substitution gives endpoint maximum
\(\ell(r-1)\binom n2+n(L-\delta_{k+1})\). Thus
\[
 M=\ell(r-1)\binom n2+nL,\qquad k_h=(r-1+h-n)\bmod r
\]
gives the unique maximizing vertex path. Any nonvertex uses another vertex
with positive weight in a convex decomposition, and has smaller volume.

In the rotating order based at \(k_h\), define the integer vector
\[
 \alpha_i^h=x^{*,h}_{k_h-i}-x^h_{k_h-i}\quad(0\le i<r),\qquad
 b=(s_{r-1},s_{r-2},\ldots,s_1,0).
 \tag{5}\label{frame}
\]
Successive coordinates decrease and the spread is at most \(\ell\).
Indeed the consecutive differences are the non-top slack coordinates and
their sum is \(\ell-a^h_{k_h}\). The top slices at \(h,h+1\) differ by
\(\ell\) in every row except \(k_h+1\), where they touch: the slacks
fix the relative differences and tightness fixes the smallest difference
at zero. Reading the last top slice gives \(b\). Thus the full physical
conditions, in addition to order and spread, are
\[
 T\alpha=(\alpha_{r-1},\alpha_0-\ell,\ldots,\alpha_{r-2}-\ell)
       \le\beta,\qquad \alpha^n\le b,
 \tag{6}\label{constraints}
\]
with at least one equality in each comparison. Here \(\alpha,\beta\) are
successive slices. Subtracting the vectors from \(x^*\) gives the inverse:
the terminal cap makes the last physical lengths nonnegative, and the
adjacent inequalities nest all earlier lengths. Contacts give tightness.

In fact \(\alpha_i\ge\beta_i\) for every \(i\). At contact zero,
\(\beta_0=\alpha_{r-1}\). At a contact \(p>0\), for \(i<p\) we have
\(\beta_i\le\beta_p+\ell=\alpha_{p-1}\le\alpha_i\), and for \(i\ge p\)
we have \(\beta_i\le\beta_p=\alpha_{p-1}-\ell\le\alpha_{r-1}\).
Write \(\mass\alpha=\sum_i\alpha_i\), including signs. Then
\(d=\sum_h\mass{\alpha^h}\ge n\mass{\alpha^n}\).

\section{The forward map and the bounded inverse}\label{bijection}
Suppose first that \(\alpha^n\ge0\). Then all vectors are nonnegative, and
\(\alpha^n_{r-1}=0\) by the cap. At a given deficit \(d\), use the same
padding length \(N=\max(d,n+1)\) for every object. Append zero vectors
through slice \(N\). The first new comparison has contact at zero because
\(\alpha^n_{r-1}=0\) and \(\alpha^n_0\le\ell\). In each slice set
\(u_h=\alpha^h_{r-1}\). Encode the vector by the multiset \(A_h\) of
\(\ell\) letters in \(\{0,\ldots,r-1\}\) determined by
\[
  \alpha_i^h=u_h+\#\{a\in A_h:a>i\}.
\]
These letters are the column heights of the diagram with row lengths
\(\alpha_i^h-u_h\), padded with zeros. Write \(m(A,B)\) for the maximum
size of a matching of distinct occurrences \(a<b\) from the two rows.
Then
\[
 m(A,B)=\min_{0\le p<r}\big(\#\{a\in A:a<p\}+\#\{b\in B:b>p\}\big).
\]
Each threshold covers every edge. In any cover, the least uncovered left
letter \(p\) forces all left letters below \(p\) and all right letters above
\(p\) into the cover. If every left letter is covered, the threshold
\(p=r-1\) costs at most the size of that cover. The bipartite matching
theorem gives the equality. Substitution in \eqref{constraints} shows that
\(u_h-u_{h+1}\) is bounded above by every threshold, with equality for a
contact. Since \(u_N=0\), the deficit becomes
\[
 d=\sum_{h=1}^N\sum_{a\in A_h}a+
           r\sum_{h=1}^{N-1}h\,m(A_h,A_{h+1}).
 \tag{7}\label{energy}
\]

\paragraph{The exchange.}
Move a singleton \(b\) to the left across a wide row \(A\) by removing
its largest letter \(a<b\), or its largest letter if there is no strict
predecessor. Output the singleton \(a\), then the wide row
\(A-\{a\}+\{b\}\). To undo the exchange, remove the smallest strict
successor of \(a\) in that row, wrapping to its smallest letter if needed.
Insert \(a\) in its place and output the removed letter as the right singleton.
An equal-width exchange leaves the contents of the two positions fixed.
These are the Nakayashiki--Yamada one-row rules \cite[Rule~3.11]{NY}
in the conventions of \cite[\S2.3]{FOY}. Their local energy
in that convention is \(-m\). Indeed, process right letters from largest to
smallest. A largest available strict predecessor can serve its partner in a
maximum matching by exchange. If there is none, no remaining right letter
has one. For a left singleton the reverse rule has the same strict-pair
count. This identifies the normalization in \cite[Definition 2.5]{FOY}.

We use a sum of transported local energies, as in \cite[Definition~3.4 and Remark~3.1]{SW}.
Start a separate scan at each position \(i\) with its row as carrier. Read
the original suffix to its right. At each row \(B\), score \(m(P,B)\) for
the current carrier \(P\), exchange the pair, and continue with the right
output. Let \(E\) be the sum of the scores of all scans. Then
\(\mathcal D=\sum_{\text{all letters}}a+rE\) is invariant under adjacent
exchanges. For a carrier starting before a swapped pair, the affine
Yang--Baxter identity gives the same exiting carrier with the same spectral
shift, hence the same sum of its two scores \cite[Proposition 2.6]{FOY}.
If \((A,B)\mapsto(B',A')\), the scans starting at the pair count the same
mutual match in either order, and enter the remaining suffix as \(A',B\)
before exchange and as \(B,A'\) afterwards, by invertibility. All later
scans are unchanged. On equal widths the carrier becomes the row just read,
so each adjacent score at \(h\) occurs \(h\) times as in \eqref{energy}.

Appending zero rows costs nothing, since \(m(P,0)=0\). If the last row
containing a positive letter is at position \(t\), then \(\mathcal D\ge t\).
Each nonzero row through \(t\) contributes at least one letter-unit. Each
zero row before \(t\) starts a zero carrier, which remains zero across
intervening zero rows and scores at its first nonzero row. These scans
are distinct.

\paragraph{The maps.}
Append \(N\) zero singletons to the \(N\) wide rows. Move the singletons,
in their original order, across the wide block to its left. At \(d\le N\)
the support bound forces every trailing wide row to be zero. Read the output
singletons as \(y_1,\ldots,y_N\), put \(y_{N+1}=0\), and output the least
nonnegative residues \(e_h=(y_h-y_{h+1})\bmod r\). Then
\[
 \sum_{h=1}^N h e_h=\sum_{h=1}^N y_h+
        r\sum_{h=1}^{N-1}h\ind_{y_h<y_{h+1}}=d.
\]
Thus the \(e_h\) are multiplicities of a partition of \(d\).

Conversely, take any partition of \(d\) with multiplicities \(e_h<r\).
Its parts are at most \(N\). Recover
\(y_h=(\sum_{a=h}^N e_a)\bmod r\), append \(N\) zero wide rows, and undo
the swaps in reverse order. All trailing singletons vanish by the same
support bound. From the \(N\) recovered wide rows set \(u_N=0\) and
\(u_h=u_{h+1}+m(A_h,A_{h+1})\). The column formula gives nonnegative
vectors of total weight \(d\), with order, spread and adjacent contacts.
The row differences recover each multiset, the terminal value \(u_N=0\)
and adjacent matches fix every baseline, and each swap is invertible.
Thus the two maps are inverse on these padded chains. A partition comes
from the ordinary sector at bound \(n\) exactly when its recovered chain
vanishes after slice \(n\) and lies below the cap at \(n\) with contact.

\begin{lemma}\label{cutoff}
If \(d\le mn\), all recovered slices after \(n\) are zero and the
retained terminal slice fits its physical cap with contact.
\end{lemma}
\begin{proof}
If a nonzero \(\beta\) follows \(\alpha\), then
\(\mass\alpha\ge m=\min(r,\ell)\). Otherwise \(\alpha_{r-1}=0\) and
\(\alpha_0<\ell\). Every noninitial coordinate of \(T\alpha\) is negative,
while its initial coordinate is zero and \(\beta_0>0\). Contact is impossible.
By coordinatewise nesting, occupation in slice \(h\) therefore costs
at least \(m(h-1)+1\). At \(d\le mn\) this excludes every \(h>n\).
Now \(T\alpha^n\le0\) gives \(u_n=0\) and nesting gives
\(\mass{\alpha^n}\le\lfloor d/n\rfloor\le m\).
For \(1\le k<r\), the integer cap \(b_{k-1}=s_{r-k}>\ell(r-k)/r\)
is positive and at least \(m-k+1\), since \(\ell(r-k)/r\ge m-k\).
If the \(k\)-th entry \(a=\alpha^n_{k-1}\) exceeds its cap, then
\(a+k\ge m+2\). For \(k\ge2\) we also have \(a\ge2\), so the first
\(k\) entries have total mass at least \(ka\ge k+a\ge m+2\).
For \(k=1\) a violation already costs \(m+1\). Both exceed the available
mass \(m\). The final coordinate is zero on both sides, giving contact.
\end{proof}

\paragraph{The first missing partitions.}
At \(d=mn+1\), if the recovered chain extends past \(n\), the occupation
bound must be an equality. Thus \(\alpha^{n+1}=\varepsilon=(1,0,\ldots,0)\),
all later slices vanish, and the first \(n\) each have mass \(m\).
By nesting those \(n\) slices are equal, with common value \(\lambda\).
At contact zero in \(\lambda\to\varepsilon\), the last entry of
\(\lambda\) is one, so its mass forces \(m=r\) and \(\lambda=(1,\ldots,1)\).
At any other contact some entry equals \(\ell\), forcing \(m=\ell\)
and \(\lambda=(\ell,0,\ldots,0)\). Conversely, each displayed constant
block followed by \(\varepsilon\) and zeros satisfies order, spread and
all adjacent contacts. In the first case the block touches itself and
\(\varepsilon\) at coordinate zero, and in the second at coordinate one.

The digit image is explicit. When \(r<\ell\), the padded wide rows
contain just one positive letter, a one in slice \(n+1\).
The output singleton word also has just one one. If it is at position
\(a\), its weight is \(1+r(a-1)\) by the scan formula for singletons
(also valid at \(a=1\)). Hence \(a=n+1\), giving \(e_n=r-1\) and
\(e_{n+1}=1\). When \(\ell<r\), there are \(\ell n+1\) ones. Their letter
sum already equals \(d\), so the singleton word has no ascents and
consists of \(\ell n+1\) initial ones. Its partition is the single part
\((\ell n+1)\). In either case the unique reverse image occupies slice
\(n+1\).

If the tail after \(n\) is empty, \(u_n=0\) and
\(\mass{\alpha^n}\le\lfloor(mn+1)/n\rfloor\). For \(n\ge2\) the mass is
at most \(m\), so the cap holds. For \(n=1\) the available mass is \(m+1\).
The estimate in Lemma~\ref{cutoff} permits a violation only at the first
coordinate. It forces
\(\alpha^1=(m+1,0,\ldots,0)\) and \(b_0=m\). The spread bound
forces \(m=r<\ell\), and \(r=b_0>\ell(r-1)/r\) forces \(\ell=r+1\).
Here \(b_0=s_{r-1}=r\) means \(c_{r-1}=1\). Conversely this vector
violates precisely that cap. Its \(\ell=r+1\) letters are ones and their
sum accounts for the full weight, so its image is the single part
\((r+1)\), distinct from \((2,1,\ldots,1)\). This exhausts all ways the
inverse can fail at \(mn+1\).

Write \(R_n^+(z)=\sum_x z^{M-W(x)}\), summing over the ordinary sector.
The bijection through deficit \(n+1\) gives
\[
 R_n^+(z)\equiv D_n(z)+z^{n+1}\pmod{z^{n+2}},\qquad
 D_n(z)=\prod_{h=1}^n(1+z^h+\cdots+z^{(r-1)h}).
 \tag{8}\label{digit}
\]
The extra object has \(e_{n+1}=1\).

\section{The signed configurations and the first gap}\label{signed}
For a signed terminal vector \(\gamma=\alpha^n\) one has \(\gamma_{r-1}<0\).
Choose a contact at \(p<r-1\), set \(j=r-1-p\), and write \(s=s_j=b_p\).
Then \(s<\ell\), since \(\gamma_{r-1}\ge s-\ell\). Order and spread give
a coordinatewise lower bound
\[
 \gamma\ge v_j=(s^{[r-j]},(s-\ell)^{[j]}),\qquad
 \mass{\gamma}\ge rs-\ell j=\delta_j>0 ,
 \tag{9}\label{contact}
\]
where brackets denote repetition. Equality of sums forces equality in
every coordinate. Conversely, for \(s_j<\ell\) this \(v_j\) fits \(b\):
the upper block is at most the corresponding caps and the lower block is
negative. It touches at \(p\), and the constant chain satisfies
\(Tv_j\le v_j\) with contact at \(p+1\).

Thus a signed contact exists precisely when \(E\ne\varnothing\).
Emptiness means \(s_1=\ell\). Otherwise the scores are distinct by
\eqref{scores}, so their minimum on \(E\) has a unique index \(j_*\).
For any contact \(j\) we have
\(\alpha^h\ge\gamma\ge v_j\) and \(d\ge n\delta_j\ge n\kappa\).
At equality every slice equals \(v_{j_*}\), and its jump of size \(\ell\)
places all the slack at \(k_h+j_*\bmod r\). This constant chain therefore
gives the unique first signed cylinder. Score one is possible only at the
index \(j\) satisfying \(rt-\ell j=1\) and only if \(s_j=t\).
Conversely \(0<t<\ell\) makes this index eligible whenever \(s_j=t\).

At \(d=n\kappa+1\), \(n\ge2\), any contact has \(j=j_*\), since
\(n(\kappa+1)>d\). Nesting forces \(\gamma=v_{j_*}\) and confines the
extra unit to the first slice. Order and spread allow only the first
lower-block entry to increase. Writing \(p=r-1-j_*\), this gives
\(\alpha^1=v_{j_*}+e_{p+1}\), with \(e_{p+1}\) a coordinate unit vector.
The comparison \(T\alpha^1\le v_{j_*}\) retains contact at \(p+1\).
For \(j_*=1\) the increment raises coordinate zero of \(T\alpha^1\)
from \(s-\ell\) to \(s-\ell+1\le s\), where \(s=s_{j_*}\).
Otherwise it raises a lower bound from \(s-2\ell\) to \(s-2\ell+1\).
The terminal slice stays signed.

When \(n=1\), a signed vector of weight two touches a cap with
\(\delta_j=1\) or \(2\). In the first case the same unique lift fits
the cap because the raised entry is at most zero. It remains signed:
if the lower block has length one, \(t-\ell=-(\ell(r-1)-1)/r\le-2\)
for the admissible integer \(t\). In the second case \(\gamma=v_j\)
by equality in \eqref{contact}. Here \(s_j<\ell\), since \(s_j=\ell\)
would give \(\delta_j=\ell(r-j)\ge3\). There is at most one index
of each score. The two vectors are distinct: the lift either has three
heights or has spread \(\ell-1\), while the equality vector has two
heights at distance \(\ell\).

\paragraph{Normalization.}
Removing \(f_h=\min_i(x_i^h-x_i^{h+1})\) common copies of each part \(h\)
leaves a tight cylinder; reinsertion adds \(rh f_h\) in volume. Hence
\(F_n(q)=\sum_{x\ {\rm tight}}q^{W(x)}/D_n(q)\). Its degree is
\(M-(r-1)\binom{n+1}{2}=d_n\) and its leading coefficient is one.
Set \(I=\ind_{\{\exists j:\delta_j=1\}}\), \(J=\ind_{\{\exists j:\delta_j=2\}}\).
By \eqref{digit} and \eqref{contact}, the reversed tight numerator, through
degree \(n+1\), equals
\[
 D_n(z)+z^{n+1}+I(z^n+z^{n+1})+\ind_{n=1}J z^2 .
\]
The denominator is reciprocal and begins \(1+z\), so the reversal of \(F_n\)
has expansion \eqref{formula}. Finally \(s_j\ge\lceil\ell j/r\rceil\)
for \(0<j<r\). The residues of \(\ell j\) modulo \(r\) permute
\(1,\ldots,r-1\), whence
\(\sum_{j=1}^{r-1}\lceil\ell j/r\rceil=(r-1)(\ell+1)/2\) and \(B\ge0\).
Also \(A\ge3\) since \(r,\ell\ge3\) are coprime. Inducting in \eqref{tri},
each \(k<n\) has degree at least the following distance below \(d_n\):
\[
 d_n-d_k-k(n-k)=(n-k)\big(An+(A-1)k+B\big)\ge An>n+1.
\]
This proves Corollary~\ref{coeff}.

\section{Two examples}\label{examples}
\paragraph{The original four-row deletion.}
The outgoing profile \((0,2,1,2)\) has canonical rotation \((2,1,2,0)\)
and scores \((\delta_0,\ldots,\delta_3)=(0,3,2,5)\). Here
\(M=15\binom n2+10n\), \(d_n=6n^2+n\), and the terminal cap in canonical
coordinates is \((5,3,2,0)\). The first signed chain occurs at \(2n\),
with difference \((3,3,-2,-2)\) at every height.
In the original row numbering \(0,1,2,3\), the top slack index at height
\(h\) is \(h-n\bmod4\). For \(d=1\) the only digit is \(e_1=1\);
the wide rows are \(00001\) followed by zero rows. The first difference
vector is \((1,0,0,0)\) and all later ones vanish. Thus the physical
operation deletes a single rightmost \(1\) in row \(1-n\bmod4\).
For example, the actual maximum row lengths are
\[
\begin{array}{c c c c c}
\toprule n & x^{*,1} & x^{*,2} & x^{*,3} & W\text{ after deletion}\\
\midrule 2 &(5,5,7,8)&(5,0,2,3)&\text{---}&35-1=34\\
3 &(10,10,12,8)&(5,5,7,8)&(5,0,2,3)&75-1=74\\ \bottomrule
\end{array}
\]
At \(n=2\) one can also see the first surviving coefficient. The extra
digit \(e_3=1\) gives singletons \(1,1,1\). The reverse exchanges give
wide rows \(00111,00000,00000\). Removing the empty last slice and
subtracting the first difference \((3,0,0,0)\) gives
\(x^1=(5,5,7,5)\), \(x^2=(5,0,2,3)\), a tight cylinder of volume \(32\).
This is the object left over after cancellation at deficit three.

\paragraph{Two profiles ending in one.}
At \(\ell=r+1\), the distinguished prefix has \(j=r-1\), \(t=r\),
so \(s_j=t\) means exactly that the last entry is one. At \((r,\ell)=(4,7)\)
one has \(j=1,t=2\). The two canonical profiles below have the same last
entry and different gaps.
\[
\begin{array}{c c c c}
\toprule c&(\delta_1,\delta_2,\delta_3)&d_n&
\text{deficit of first surviving coefficient below the top }(n\ge2)\\ \midrule
(2,2,2,1)&(1,2,3)&9n^2&n\\
(3,1,2,1)&(5,2,3)&9n^2+n&n+1\\ \bottomrule
\end{array}
\]
For the first cap \(b=(6,4,2,0)\), contact at \(p=2\) gives
\((2,2,2,-5)\). Its constant chain is the sole signed configuration at
\(d=n\); its lift changes the first slice to \((2,2,2,-4)\).
For the second cap \(b=(6,4,3,0)\), even the cheapest signed contact
costs two per slice. At \(n=2\) the first profile gives terms
\(q^{36}+q^{34}+q^{33}\) at the top, as in Warnaar's printed polynomial
\cite[pp.~43--44]{W}. Its \(q^{35}\) coefficient is zero. This profile is also a rotation of
the balanced profile \((1,2,2,2)\) in \cite[(1.16)]{HLO}. Their
positive finite sum already applies to every \(Q_{n,c}\) in this orbit
\cite[Corollary~1.3, (1.19)]{HLO}.

At \(n=2\), both level-seven profiles have the same first failure of the
ordinary-sector inverse. Its range is complete through \(d=8\).
At \(d=9\) the partition \((3,2,2,2)\) reconstructs the slices
\((1,1,1,1),(1,1,1,1),(1,0,0,0)\). The third slice is occupied.
For a small example of the additional cap obstruction, take
\((r,\ell)=(3,4)\), \(c=(2,1,1)\), \(n=1\). At \(d=4\) the partition
\((2,1,1)\) needs two slices, while \((4)\) gives a single slice
\((4,0,0)\) exceeding the cap \((3,2,0)\). The other permissible
partitions, \((3,1)\) and \((2,2)\), both fit.

\section{An independent calculation from the bounded determinant}
Use the Foda--Welsh specialization of the Gessel--Krattenthaler determinant
\cite{GK,FW}, in the form of \cite[Theorem 2, (14)]{FW}, with offsets \(\mu_t=s_{r-1}-s_{t-1}\),
\(1\le t\le r\). The successive differences are \(c_0,\ldots,c_{r-2}\)
and the wrap is \(\ell+\mu_r-\mu_1=c_{r-1}\).
Expand along permutations \(t=\sigma(s)\), put \(m_t=n+rk_s+s-t\)
with \(\sum_s k_s=0\), and drop terms with negative Pochhammer index.
The remaining terms are indexed by \(m_t\ge0\) with sum \(rn\) and
distinct residues \(m_t-n+t\bmod r\), which recover the rows \(s\), the
permutation and its integers \(k_s\). Each term of \(F_n\) is
\(\operatorname{sgn}(\sigma)q^H(q;q)_n/\prod_t(q;q)_{m_t}\) with
\[
 H=\sum_{t=\sigma(s)}\left((\mu_t-t)(m_t-n)+
       \frac{r+\ell}{2r}\big((m_t-n+t)^2-s^2\big)\right).
\]
Let \(E(m)=H+\binom{n+1}{2}-\sum_t\binom{m_t+1}{2}\) and
\(E_*=E(rne_1)\). In \(E(m)\) the coefficient of \(m_t^2\) is \(\ell/(2r)\);
the linear coefficient is \(\mu_t+\ell t/r\) up to a constant independent
of \(t\). Using \(r(\mu_t-\mu_1)+\ell(t-1)=-\delta_{t-1}\) gives
\[
 E_*-E(m)=\frac{\ell}{2r}\big((rn)^2-\sum_t m_t^2\big)
                 +\frac1r\sum_t m_t\delta_{t-1}.
\]
Substitution of \(rne_1\) gives \(E_*=An^2+Bn\), using
\(\sum_t\mu_t=r s_{r-1}-L\). All vertices \(rne_t\) occur. Their common
permutation is a cyclic shift by \(n\), of sign \((-1)^{(r-1)n}\),
which cancels the Pochhammer reversal sign. Put \(Q^*(z)=z^{d_n}Q_{n,c}(1/z)\).
The reversed vertex sum is
\[
 V(z)=\frac{\sum_{j=0}^{r-1}z^{n\delta_j}}{\prod_{a=n+1}^{rn}(1-z^a)} .
\]
Every nonvertex integer composition has \(\sum m_t^2\le(rn-1)^2+1\)
and hence loss at least \(\ell(n-1/r)>n+1\). The last inequality holds
at \(n=1\) since \(r,\ell\ge3\) are coprime, and its margin increases with
\(n\). Expanding \(V\) and applying the degree separation in
Section~\ref{signed} gives \eqref{formula} independently.

\begin{proposition}\label{vertex}
Put
\[
 H_i=n(\ell+\delta_i)-c_{i-1}\quad(0\le i<r),\qquad
 K=\min_{0\le i<r} H_i ,
\]
using cyclic indices. Then \(Q^*(z)\equiv V(z)\pmod{z^K}\), and
\[
 [z^K](Q^*-V)=-\#\{i:H_i=K\}-\ind_{\{n=1,\ K=d_1\}} .
\]
\end{proposition}
\begin{proof}
Set \(S=rn\) and index columns from zero. Start at a vertex \(Se_i\)
and move one unit to column \(j\ne i\). The two affected row residues
change by \(-1,+1\), respectively. To remain distinct they must exchange
their former residues, which occurs exactly when \(j=i-1\bmod r\).
Each such composition is therefore admissible and has the opposite
permutation sign. Its reversed Pochhammer factors have the same leading
sign as at a vertex, since \(\sum m_t=S\). Its leading coefficient is
\(-1\) and the loss formula gives \(H_i\), using
\(\delta_i-\delta_{i-1}=r c_{i-1}-\ell\).

Every other nonvertex has \(\sum_{a<b}m_am_b\ge S\). With two positive
parts both are at least two, so their product is at least their sum.
With three or more, sum \(2m_am_b\ge m_a+m_b\) over pairs in the support.
The score term is strictly positive since only \(\delta_0=0\).
Thus these losses exceed \(\ell n\ge H_0\).

In \eqref{tri} the loss of the term indexed by \(k<n\) is
\(An^2+Bn-(n+B)k-(A-1)k^2\), strictly decreasing in \(k\).
The smallest loss \(T=(2A-1)n-A+1+B\) therefore comes from \(k=n-1\),
with leading coefficient one before subtraction. To compare it with \(K\),
use \(s_{r-1}=\ell-c_{r-1}\) to write
\[
 T-H_0=\sum_{j=1}^{r-2}s_j-(r-1)+(n-1)\bigl((r-2)\ell-r\bigr).
\]
Each prefix in the sum is at least one, and \(s_{r-2}>\ell(r-2)/r\ge1\).
The slope \((r-2)\ell-r\) is positive for coprime \(r,\ell\ge3\).
Hence \(T\ge H_0\ge K\), strictly in the first inequality for \(n\ge2\).
Only at \(n=1\), when \(T=d_1\), can the subtraction contribute at \(K\).
\end{proof}

For both level-seven profiles in Section~\ref{examples} at \(n=2\)
the physical inverse first misses a partition at deficit nine, while
the vertex sum first differs from \(Q^*\) at \(K=13\).
The coefficients of \(Q^*\) and \(V\) there are respectively \(16,17\)
for \((2,2,2,1)\) and \(12,13\) for \((3,1,2,1)\).

\par\medskip
\noindent Independent researcher\\
\href{mailto:v@deltagray.com}{\nolinkurl{v@deltagray.com}}
\end{document}